\documentclass{article}
\usepackage{amssymb}
\usepackage{amsmath}

\newtheorem{theorem}{Theorem}
\newtheorem{lemma}{Lemma}

\begin{document}

\addvspace{2in}

\title{Global Existence and Increased Spatial Decay for the Radial Vlasov-Poisson System with Steady Spatial Asymptotics}
\author{Stephen Pankavich \\
Department of Mathematics \\
Indiana University \\
Bloomington, IN 47401 \\
sdp@indiana.edu} \date{} \maketitle

\begin{abstract}
A collisionless plasma is modelled by the Vlasov-Poisson system in
three space dimensions.  A fixed background of positive charge,
which is independent of time and space, is assumed.  The situation
in which mobile negative ions balance the positive charge as $\vert
x \vert \rightarrow \infty$ is considered.  Hence the total positive
charge, total negative charge, and total energy are infinite. Smooth
solutions with appropriate asymptotic behavior for large $\vert x
\vert$, which were previously shown to exist locally in time, are
continued globally for spherically symmetric data. This is done by
showing that the charge density decays at least as fast as $\vert x
\vert^{-4}$. Finally, an increased decay rate of $\vert x
\vert^{-6}$ is shown in the general case without the assumption of
spherical symmetry.
\end{abstract}

\vspace{0.1in}

\begin{center}
INTRODUCTION \\
\end{center}

Let $F : \mathbb{R}^3 \rightarrow [0,\infty)$, $f_0 : \mathbb{R}^3
\times \mathbb{R}^3 \rightarrow [0,\infty)$, and $A : [0,\infty)
\times \mathbb{R}^3 \rightarrow \mathbb{R}^3$ be given. We seek a
solution, $f : [0,\infty) \times \mathbb{R}^3 \times \mathbb{R}^3
\rightarrow [0,\infty)$ satisfying

\begin{equation}
\label{one} \left. \begin{array}{ccc}
& & \partial_t f + v \cdot \nabla_x f - (E + A) \cdot \nabla_v f = 0,\\
\\
& & \rho(t,x) = \int (F(v) - f(t,x,v)) dv,\\
\\
& & E(t,x) = \int \rho(t,y) \frac{x - y}{\vert x - y \vert^3} dy,\\
\\
& & f(0,x,v) = f_0(x,v). \end{array} \right \}
\end{equation}
Here $F$ describes a number density of positive ions which form a
fixed background, and $f$ denotes the density of mobile negative
ions in phase space.  Notice that if $f_0(x,v) = F(v)$ and $A =
0$, then $f(t,x,v) = F(v)$ is a steady solution.  Thus, we seek
solutions for which $f(t,x,v) \rightarrow F(v)$ as $\vert x \vert
\rightarrow \infty$. Precise conditions which ensure local
existence were given in \cite{VPSSA}.  It is important to notice
that (\ref{one}) is a representative problem, and that problems
concerning multiple species of ions can be treated in a similar manner.\\

The paper will be divided into two main results. The first will be
devoted to showing global existence of a smooth solution to
(\ref{one}) in the case of spherically symmetric data. This serves
to continue the local existence result of \cite{VPSSA}. The second,
then, is devoted to achieving an increased rate of spatial decay of
the charge density without the spherical symmetry assumption.  We
will assume throughout that $F$ has compact velocity support, so
that decay in $v$ of the background is not an issue. Thus, the main
difficulty in showing global existence arises in showing $\rho$
decays rapidly enough in $\vert x \vert$. To see this, consider the
following heuristic argument, discussed in \cite{VPSSA}.  Let $r =
\vert x \vert$. Then, we typically expect the electric field, $E$,
to decay like $r^{-2}$ for large $r$. If we let $g = F - f$, then

\begin{equation}
\label{vlasov}
\partial_t g + v \cdot \nabla_x g - (E + A) \cdot \nabla_v g
= -(E + A) \cdot \nabla_v F.
\end{equation}
Viewing $(E + A) \cdot \nabla_v F$ as a source term for $g$, we can
only conclude $r^{-2}$ decay for both $g$ and $\rho$. But, if $\rho$
really decayed like $r^{-2}$, the integral for $E$ could certainly
fail to decay like $r^{-2}$. In fact, we cannot show that $g$ decays
faster than $r^{-2}$, but due to cancellation in the $v$ integral,
$\rho$ must decay faster.\\

The Vlasov-Poisson system has been studied extensively in the case
when $F(v) = 0$ and solutions decay to zero as $\vert x \vert
\rightarrow \infty$.  Smooth solutions were shown to exist globally
in time in \cite{P} and independently in \cite{LP}. Important
results prior to global existence appear in \cite{B}, \cite{Horst},
\cite{Horst82}, and \cite{K}. Also, the method used by \cite{P} is
refined in \cite{Horst93} and \cite{S}. Global existence for the
Vlasov-Poisson system in two dimensions was established in \cite{OU}
and \cite{W}. A complete discussion of the literature concerning
Vlasov-Poisson may be found in \cite{G}. We also mention \cite{BR}
since the problem treated in said paper is periodic in space, and
thus the solution does not decay for large $\vert x \vert$.  The
works cited above make extensive use of the laws of conservation of
charge and energy.  However, in the problem considered here (and
those of \cite{C}, \cite{J}, \cite{SSAVP}, and \cite{VPSSA}), the
charge and energy are infinite, and it is less clear how to use the
conservation laws. Therefore, the use of conservation laws will be
an important issue, and we will utilize a lemma (stated here as
Lemma $2$) from \cite{SSAVP} to deal with it properly.

\vspace{.25in}

\section*{Preliminaries}

\vspace{.15in}

We will use notation which follows \cite{VPSSA}. For $q > 7 +
\sqrt{33}$, let
$$ p  = 4 - \frac{8}{q} $$
and denote
$$R(x) = R(\vert x \vert) = (1 + \vert x \vert^2)^\frac{1}{2}.$$
We will use the norms

$$\Vert g \Vert_\infty = \sup_{z \in
\mathbb{R}^n} \vert g(z) \vert $$

$$ \Vert \rho \Vert_p = \Vert \rho R^p(x)
\Vert_\infty,$$

$$ \Vert g \Vert_q = \Vert g (1 + \vert x \vert^2  +\vert v \vert^q)
\Vert_{L^\infty(\mathbb{R}^6)}$$ and
$$ \vert \vert \vert g \vert \vert \vert = \Vert g \Vert_q + \Vert \nabla g \Vert_q +
\left \Vert \int g \ dv \right \Vert_p,$$ but never use $L^p$ or
$L^q$ for finite $p$ and $q$.  We will write, for example, $\Vert
g(t) \Vert_q$ for the $\Vert \cdot \Vert_q$ norm of $(x,v) \mapsto
g(t,x,v)$.  Notice that we may take $q$ arbitrarily large since we
take the initial
data to be of compact support.\\

Following \cite{SSAVP} and \cite{VPSSA}, we assume the following
conditions hold for some $C > 0$ and all $t \geq 0$, $x \in
\mathbb{R}^3$, and $v \in \mathbb{R}^3$, unless otherwise stated :

\begin{enumerate}
\renewcommand{\labelenumi}{(\Roman{enumi})}
\item $F(v) = F_R(\vert v \vert)$ is nonnegative and $C^2$ with,

\begin{equation}
F^{''}_R(0) < 0.
\end{equation}

In addition, there is $W \in (0,\infty)$ such that

\begin{equation}
\left . \begin{array}{ccc}
F_R^\prime(u) < 0 \ \ \ \ \mathrm{for} \ u \in (0,W)\\
\\
F_R(u) = 0 \ \ \ \ \mathrm{for} \ u \geq W.
\end{array} \right \}
\end{equation}

\item $f_0$ is $C^1$ with compact $v$-support, nonnegative, and satisfies the condition of
spherical symmetry,

\begin{equation}
\label{f0} f_0(x,v) = f_R (\vert x \vert, \vert v \vert, x \cdot
v)
\end{equation}

which is equivalent to

$$ f_0(x,v) = f_0(Ux, Uv)$$

for every rotation $U$.

Also, there is $N > 0$ such that for $\vert x \vert > N$, we have
$$ f_0(x,v) = F(v).$$

\item $A$ is $C^1$ and

$$\vert A(t,x) \vert + \vert \partial_x A(t,x) \vert \leq C R^{-2}(x)$$

$$ \vert \nabla_x \cdot A(t,x) \vert \leq C R^{-4}(x). $$
Furthermore, $A$ satisfies the condition of
spherical symmetry,

\begin{equation}
\label{A}
A(t,x) = \alpha(t, \vert x \vert) \frac{x}{\vert x
\vert^2}
\end{equation}
which is equivalent to

$$ A(t,x) = U^T A(t, Ux)$$
for every rotation $U$.\\
Finally, we assume that there is a continuous function $a : [0,T]
\rightarrow \mathbb{R}$ such that for $\vert x \vert > N$
$$\left \vert \alpha(t,\vert x \vert) - \frac{a(t)}{\vert x \vert} \right \vert \leq R^{2-p}(x).$$
\end{enumerate}

It should be noted that the assumptions (\ref{f0}) and (\ref{A}),
imply the spherical symmetry of $f(t,x,v)$, $g(t,x,v)$, and $E(t,x)$
for all $t \in [0,\infty)$, $x \in \mathbb{R}^3$, and $v \in
\mathbb{R}^3$. Thus, where necessary we will write $g(t,x,v) = g(t,
\vert x \vert, \vert v \vert, x \cdot v)$ when using the spherical
symmetry of $g$. Furthermore, the spherical symmetry of $E$ and $f$
will be instrumental in the first results of the paper, while the
last theorem is dedicated to eliminating the symmetry assumptions to
conclude similar decay rates. We wish to prove the following :

\begin{theorem}
Assuming conditions $(I)$, $(II)$, and $(III)$ hold, there exists
$f \in \mathcal{C}^1([0,\infty) \times \mathbb{R}^3 \times
\mathbb{R}^3)$ that satisfies (\ref{one}) with $\vert \vert \vert
(F - f)(t) \vert \vert \vert$ bounded on $t \in [0,T]$, for every
$T > 0$. Moreover, $f$ is unique.
\end{theorem}

\vspace{.25in}

In \cite{VPSSA}, both local existence and a criteria for
continuation of a bounded, unique solution of (\ref{one}) are
shown. Thus, to prove Theorem $1$ we will need to establish
the continuation criteria for all $T > 0$.  Specifically, Theorems
$2$ and $3$ of \cite{VPSSA}, when combined, state :

\begin{theorem}
Assume $q > 7 + \sqrt{33}$ and conditions $(I)$, $(II)$ and
$(III)$ hold, without (\ref{f0}) and (\ref{A}). Let $f$ be a
$\mathcal{C}^1$ solution of (\ref{one}) on $[0,T] \times
\mathbb{R}^3 \times \mathbb{R}^3$ with $T > 0$. If $\Vert \int (F
- f)(t) \ dv \Vert_p$ is bounded on $[0,T]$, then we may uniquely
extend the solution to $[0,T + \delta]$ for some $\delta
> 0$ with $\vert \vert \vert (F - f)(t) \vert \vert \vert$ bounded
on $[0,T + \delta]$.
\end{theorem}

\vspace{.25in}

Therefore, to prove Theorem $1$, we will find a solution to
(\ref{one}) on $[0,T]$ for some $T > 0$ using the local existence
theorem (again, see \cite{VPSSA}), that has been uniquely extended
using Theorem $2$. Since this may be done as long as the
$p$-norm stays bounded, we will only need to prove the following
lemma to establish Theorem $1$.

\begin{lemma}
Assume $q > 7 + \sqrt{33}$ and conditions $(I), (II)$ and $(III)$
hold. Let $f$ be a $\mathcal{C}^1$ solution of (\ref{one}) on
$[0,T] \times \mathbb{R}^3 \times \mathbb{R}^3$. Then,

$$ \left \Vert \int (F - f)(t) \ dv \right \Vert_p \leq C $$
for all $t \in [0,T]$, where $C$ is determined by $F$, $A$, $f_0$,
and $T$.
\end{lemma}

\vspace{.3in}

Although it is not explicitly stated in \cite{VPSSA}, the proof
presented there shows that $\delta$ in Theorem $2$ is bounded away
from zero as long as $\Vert \int (F - f)(t) \ dv \Vert_p$ is
bounded.  Using this observation, Theorem $1$ follows from Lemma
$1$.\\

Once Lemma $1$ has been established, we will show increased decay of
the charge density, $\rho$, under slightly modified assumptions.
Instead of $(III)$, we will assume the following conditions for some
$C > 0$ and all $t \geq 0$, $x \in \mathbb{R}^3$, and $v \in
\mathbb{R}^3$:

\begin{enumerate}
\item[(IV)] $A$ is $C^1$  with

\begin{equation}
\vert A(t,x) \vert \leq C R^{-2}(x),
\end{equation}

\begin{equation}
\vert \partial_x A(t,x) \vert \leq C R^{-3}(x),
\end{equation}

and

\begin{equation}
\vert \nabla_x \cdot A(t,x) \vert = 0.
\end{equation}

\end{enumerate}

Notice that condition $(IV)$ does not involve spherical symmetry of
$A$, and thus the result which follows requires no such assumption.

\begin{theorem}
Assume conditions $(I)$ and $(II)$ hold, without $(\ref{f0})$, and
condition $(IV)$ holds. Let $T > 0$ and $f$ be the \ $\mathcal{C}^1$
solution of (\ref{one}) on $[0,T) \times \mathbb{R}^3 \times
\mathbb{R}^3$ with
$$ \vert \vert \vert (F - f)(t) \vert \vert \vert < \infty $$
for every $t \in [0,T)$.
Then, we have

$$ \Vert \rho(t) \Vert_6 \leq C_{p, t}  $$
for any $t \in [0,T)$, where $C_{p, t}$ depends upon
$$\sup_{\tau \in [0,t]} \Vert \rho(\tau) \Vert_p.$$

\end{theorem}
\vspace{0.25in}

In Section $1$ we will prove Lemma $1$, and thus Theorem $1$. Then,
Section $2$ will contain the proof of Theorem $3$.  We will denote
by ``$C$'' a generic constant which changes from line to line and
may depend upon $f_0$, $A$, $F$, or $T$, but not on $t$, $x$, or
$v$. When it is necessary to refer to a generic constant which
depends upon other variables, we will use variable subscripts to
distinguish them. For example, we will use $C_{p, t}$ quite
frequently in Section $3$, and it will always denote dependence upon
$\Vert \rho(t) \Vert_p$ and $t$.  When it is necessary to refer to a
specific constant, we will use numeric superscripts to distinguish
them.  For example, $C^{(1)}$ will always refer to the same
constant.

\vspace{.25in}

\section*{Section 1}

Define the characteristics $X(s,t,x,v)$ and $V(s,t,x,v)$ by

\begin{equation}
\label{char} \left. \begin{array}{ccc} & &
\frac{\partial X}{\partial s} (s,t,x,v) = V(s,t,x,v) \\
\\
& & \frac{\partial V}{\partial s} (s,t,x,v) = - \left( E(s,X(s,t,x,v)) + A(s, X(s,t,x,v)) \right) \\
\\
& & X(t,t,x,v) = x \\
\\
& & V(t,t,x,v) = v. \end{array} \right \}
\end{equation}
Then, we have
$$ \frac{d}{ds} f(s,X(s,t,x,v),V(s,t,x,v)) =
\partial_t f + V \cdot \nabla_x f - ( E + A ) \cdot \nabla_v f = 0.$$
Therefore, $f$ is constant along characteristics, and

\begin{equation}
\label{fchar} f(t,x,v) = f(0, X(0,t,x,v), V(0,t,x,v)) =
f_0(X(0,t,x,v),V(0,t,x,v)).
\end{equation}
Thus, we find by $(II)$ that $f$ is nonnegative and $ \sup_{x,v} f =
\Vert f_0 \Vert_\infty < \infty$.  Unless necessary, we will omit
writing the dependence of $X(s)$ and $V(s)$ on
$t$, $x$, and $v$ for the remainder of the paper.\\

In order to bound the electric field, and thus the velocity support,
we must use Lemma $3$ and Theorem $4$ from \cite{SSAVP}. In
particular, we state the following lemma without proof.

\begin{lemma}
Assuming conditions $(I)$ and $(III)$ hold, without (\ref{A}),
there exists $k : [0,T] \times \mathbb{R}^3 \rightarrow
[0,\infty)$ such that
$$\vert \rho(t,x) \vert \leq C(k(t,x)^\frac{3}{5} + k(t,x)^\frac{1}{2}) $$
and
$$ \int k(t,x) dx \leq C. $$
\end{lemma}

\vspace{0.25in}

\noindent Furthermore, notice that, due to the radial symmetry, we
may write
\begin{equation}
\label{e}
E(t,x) = \frac{m(t,\vert x \vert)}{\vert x \vert^2}
\frac{x}{\vert x \vert}
\end{equation}
where the enclosed charge is given by
$$ m(t,r) := \int_{\vert y \vert < r} \rho(t,y) dy.$$
Lemma $2$ will be exactly what is needed to bound the electric
field, $E$.  We have for any $x \in \mathbb{R}$

\begin{eqnarray*}
\int_{\vert y \vert \leq \vert x \vert} k(t,y)^\frac{1}{2} dy & \leq
& \left (\int_{\vert y \vert \leq \vert x \vert} k(t,y) dy
\right)^\frac{1}{2} \left (\int_{\vert y
\vert \leq \vert x \vert} dy \right)^\frac{1}{2} \\
& \leq & \left (\int_{\vert y \vert \leq \vert x \vert} k(t,y)
dy \right)^\frac{1}{2} \left(C \vert x \vert^3 \right)^\frac{1}{2} \\
& \leq & C \vert x \vert^\frac{3}{2}
\end{eqnarray*}
and

\begin{eqnarray*}
\int_{\vert y \vert \leq \vert x \vert} k(t,y)^\frac{3}{5} dy & \leq
& \left (\int_{\vert y \vert \leq \vert x \vert} k(t,y)
dy \right)^\frac{3}{5} \left (C \vert x \vert^3 \right)^\frac{2}{5} \\
& \leq & C \vert x \vert^\frac{6}{5}.
\end{eqnarray*}
Now, we use the bounds on $k$ to estimate the enclosed charge.

\begin{eqnarray*}
\vert m(t, \vert x \vert ) \vert & \leq &  \int_{\vert y \vert
\leq \vert x \vert} \vert \rho(t,y)\vert dy \\
& \leq & C\int_{\vert y \vert \leq \vert x
\vert}(k(t,y)^\frac{3}{5} + k(t,y)^\frac{1}{2}) \ dy \\
& \leq & C (\vert x \vert^\frac{3}{2} + \vert x \vert^\frac{6}{5}). \\
\end{eqnarray*}
Thus,
\begin{eqnarray*}
\vert E(t,x) \vert & \leq & C \vert m(t,\vert x \vert) \vert \vert
x \vert^{-2} \\
& \leq & C \vert x \vert^{-2} (\vert x \vert^\frac{6}{5} + \vert
x \vert^\frac{3}{2}) \\
& \leq & C (\vert x \vert^{-\frac{4}{5}} + \vert x
\vert^{-\frac{1}{2}}).
\end{eqnarray*}
So,
$$ \vert E(t,x) \vert \leq C \mathcal{G}(\vert x \vert),$$
and combining this with $(III)$,
$$ \vert E(t,x) + A(t,x) \vert \leq C \mathcal{G}(\vert x \vert)$$
where $$\mathcal{G}(r) := \left \{ \begin{array}{l} r^{-\frac{4}{5}}, \ \ r \leq 1 \\
r^{-\frac{1}{2}}, \ \ r \geq 1.
\end{array} \right. $$

Next, we use methods from \cite{Horst82} to estimate the velocity
support. Assume
$$f(t,x,v) = f(0,X(0,t,x,v), V(0,t,x,v)) \neq 0$$
so that, using the compact support in $v$ of $f_0$, we have
$$ \vert \dot{X}(0) \vert = \vert V(0) \vert \leq C. $$
For any $\epsilon > 0 $, define $a := 2 + \epsilon$, $b :=
\frac{a}{a-1}$, and
$$ B := (2\int_1^\infty \mathcal{G}(\vert x \vert)^a dx ) ^\frac{1}{a}.$$
Let $t_1, t_2 \in [0,T]$ with $t_1 < t_2$, and assume
$\dot{X}_i \geq 0$ on $ [t_1,t_2]$.\\
We claim
\begin{equation}
\label{Claim} \int_{X_i(t_1)}^{X_i(t_2)} \mathcal{G}(\vert x \vert)
\ dx \leq B (X_i(t_2) - X_i(t_1))^\frac{1}{b} + 10.
\end{equation}
To show (\ref{Claim}), let
$$ \mathcal{G}_1 (x) := \left \{ \begin{array}{lll} \vert x \vert^{-\frac{4}{5}}
& \ & \mathrm{if} \ \ \vert x \vert\leq 1 \\ 0 & \ & \mathrm{if} \
\ \vert x \vert > 1 \end{array} \right. $$ and
$$ \mathcal{G}_2 (x) := \left \{ \begin{array}{lll} 0 & \ & \mathrm{if} \ \
\vert x \vert \leq 1 \\ \vert x \vert^{-\frac{1}{2}} & \ &
\mathrm{if} \ \ \vert x \vert > 1.
\end{array} \right. $$
Then,
\begin{eqnarray*}
\int_{X_i(t_1)}^{X_i(t_2)} \mathcal{G}(\vert x \vert) \ dx & = &
\int_{X_i(t_1)}^{X_i(t_2)} \mathcal{G}_1 (x) \ dx +
\int_{X_i(t_1)}^{X_i(t_2)} \mathcal{G}_2 (x) \ dx \\
& \leq & \int_{-1}^1 \mathcal{G}_1 (x) \ dx +
\left(\int_{X_i(t_1)}^{X_i(t_2)} (\mathcal{G}_2 (x))^a \
dx \right)^\frac{1}{a} \left (X_i(t_1) - X_i(t_2) \right )^\frac{1}{b} \\
& \leq & 10 + \left (2 \int_1^\infty (\mathcal{G}_2 (x) \right)^a \
dx)^\frac{1}{a} \left (X_i(t_1) - X_i(t_2) \right )^\frac{1}{b} \\
& \leq & 10 + B (X_i(t_1) - X_i(t_2))^\frac{1}{b}.
\end{eqnarray*}
This establishes (\ref{Claim}).  Now, using this result and
following \cite{Horst82}, we have

\begin{eqnarray*}
\vert \dot{X}_i(t_2)^2 - \dot{X}_i(t_1)^2 \vert & = & \left \vert
2 \int_{t_1}^{t_2} \dot{X}_i(s) \ddot{X}_i(s) \ ds \right \vert \\
& \leq & C \int_{t_1}^{t_2} \mathcal{G}(\vert X(s) \vert)
\dot{X}_i(s) \ ds \\
& \leq & C \int_{t_1}^{t_2} \mathcal{G}(\vert X_i(s) \vert) \dot{X_i}(s) \ ds \\
& =: & C^{(0)} \int_{X_i(t_1)}^{X_i(t_2)} \mathcal{G}(\vert x \vert) \ dx \\
& \leq & C^{(0)}B(X_i(t_2) - X_i(t_1))^\frac{1}{b} + 10C^{(0)} \\
& \leq & C^{(0)}B ((\sup_{\tau \in [t_1, t_2]} \vert
\dot{X}_i(\tau) \vert) (t_2 - t_1))^\frac{1}{b} + 10C^{(0)}.
\end{eqnarray*}
Define
$$W := \sup_{s \in [0,T]} \vert \dot{X}_i(s) \vert.$$
Then,
\begin{equation}
\label{starV} \vert \dot{X}_i(t_2)^2 - \dot{X}_i(t_1)^2 \vert \leq
C^{(0)}B(Wt_2)^\frac{1}{b} + 10C^{(0)}.
\end{equation}
Note that this holds if $\dot{X}_i \leq 0$ on $[t_1,t_2]$, as well.
Let us consider $t \in [0,T]$ and $\dot{X}_i(t) > 0$. Define
$$\bar{t} := \inf \{ \tau \geq 0 : \dot{X}_i(s) \geq 0, \ \forall s
\in [\tau, t]\}.$$ If $\bar{t} = 0 $, then by (\ref{starV})

\begin{eqnarray*}
\dot{X}_i^2(t) & \leq & \dot{X}_i^2(0) + C^{(0)}B(Wt)^\frac{1}{b} + 10C^{(0)} \\
& \leq & \dot{X}_i^2(0) + C^{(0)}B(WT)^\frac{1}{b} + 10C^{(0)}.
\end{eqnarray*}
If $\bar{t} > 0 $, then $\dot{X}_i(\bar{t}) = 0$ and by
(\ref{starV}), we have

\begin{eqnarray*}
\dot{X}_i^2(t) & \leq & \dot{X}_i^2(\bar{t}) + C^{(0)}B(Wt)^\frac{1}{b} + 10C^{(0)} \\
& = & C^{(0)}B(WT)^\frac{1}{b} + 10C^{(0)} \\
& \leq & \vert \dot{X}_i^2(0) \vert + C^{(0)}B(WT)^\frac{1}{b} +
10C^{(0)}.
\end{eqnarray*}
We may repeat this process for $\dot{X}_i(t) < 0$, so summing over
$i$ yields
$$ \vert \dot{X}(t) \vert^2 \leq \vert \dot{X}(0) \vert^2 +
3C^{(0)}B(WT)^\frac{1}{b} + 30C^{(0)}. $$ Since the right-hand side
is independent of $t$, we take the supremum and find,
$$ W^2 \leq 3C^{(0)}BT^\frac{1}{b} W^\frac{1}{b} + M$$
where $M := \vert \dot{X}(0) \vert^2 + 30C^{(0)}$.
If $3C^{(0)}BT^\frac{1}{b}W^\frac{1}{b} \leq M$, then
$$ W \leq \sqrt{2M} =: \mathcal{C}_1(T). $$
If $M \leq 3C^{(0)}BT^\frac{1}{b}W^\frac{1}{b}$, then
$$ W \leq (6C^{(0)}B)^\frac{b}{2b-1} T^\frac{1}{2b-1} := \mathcal{C}_2(T). $$
So, combining the inequalities
$$ W \leq \max\{\mathcal{C}_1(T), \mathcal{C}_2(T)\}. $$
Thus, all characteristics $V(s,t,x,v)$, along which $f$ is non-zero,
are bounded for any $s \in [0,T]$, including $V(t,t,x,v) = v$.
 Define
$$Q_t := \sup \{\vert v \vert : \exists x \in \mathbb{R}^3, \ \tau
\in [0,t] \mbox{such that} f(\tau, x, v) \neq 0 \}.$$  Notice that
$Q$ is an increasing function of time, so that we may write for
every $ s \in [0,T]$
$$\vert V(s,t,x,v) \vert \leq Q_s \leq Q_T.$$
Since the momentum is bounded, bounds on position follow.  Note that

\begin{eqnarray*}
\vert X(s) - x \vert & \leq & \int_s^t \vert V(\tau) \vert d\tau \\
& \leq & Q_TT.
\end{eqnarray*}
So, for $\vert x \vert \geq 2Q_TT$,

\begin{equation}
\label{xchar1}
\vert X(s) \vert \geq \vert x \vert - Q_T T  \geq
\frac{1}{2} \vert x \vert
\end{equation}
and

\begin{equation}
\label{xchar2}
\vert X(s) \vert \leq \vert x \vert + Q_T T \leq
\frac{3}{2} \vert x \vert.
\end{equation}
Hence, we have control on $X$ characteristics.\\

Now, to bound the charge, we must first estimate the corresponding
density. We use (\ref{vlasov}) to find
$$ \frac{d}{ds} ( g(s, X(s), V(s))) = - (E(s, X(s)) + A(s, X(s))) \cdot \nabla_v F(V(s)). $$
Thus,
\begin{equation}
\label{FI} g(t,x,v) = g(0, X(0), V(0)) - \int_0^t (E(\tau,
X(\tau)) + A(\tau, X(\tau))) \cdot \nabla_v F(V(\tau)) d\tau
\end{equation}
and by $(III)$, (\ref{e}), and (\ref{xchar1}),

\begin{equation}
\label{g}
\begin{array}{rcl}
\vert g(t,x,v) \vert & \leq & \vert
g_0(X(0, t, x, v), V(0, t, x, v)) \vert \\
\\
&  \ & \ + \int_0^t (\vert E(\tau, X(\tau, t, x, v))\vert + \vert
A(\tau, X(\tau, t, x, v)) \vert ) \ \vert
\nabla_v F(V(\tau, t, x, v)) \vert \ d\tau \\
\\
& \leq & C \vert X(0, t, x, v) \vert^{-2} + C \int_0^t \vert
X(\tau, t, x, v) \vert^{-2} \ \Vert \nabla_v F \Vert_\infty \ d\tau \\
\\
& \ & \ + C\int_0^t \vert m(\tau, X(\tau, t, x, v)) \vert \ \vert
X(\tau, t, x, v) \vert^{-2} \Vert
\nabla_v F \Vert_\infty \ d\tau \\
\\
& \leq & C \vert x \vert^{-2} \left (1 + \int_0^t \vert m(\tau,
X(\tau, t, x, v)) \vert \ d\tau \right )
\end{array}
\end{equation}
for $\vert x \vert \geq 2Q_T T$.\\

To show that $g$ decays like $r^{-2}$, we must bound the enclosed
charge $m$. For $\vert x \vert \leq 2Q_TT$,
\begin{eqnarray*}
\vert m(t,\vert x \vert) \vert & \leq & \int_{\vert y \vert \leq 2Q_TT} \vert \rho(t,y) \vert \ dy \\
& \leq & (2Q_TT)^3 \int_{\vert v \vert \leq Q_T} \Vert F - f_0
\Vert_\infty \ dv \\
& \leq & C
\end{eqnarray*}
Now, let $r = \vert x \vert$ and define $\mathcal{M}(t) := \sup_r
\vert m(t,r) \vert$. We know by (\ref{vlasov}), (\ref{g}), and the
Divergence Theorem,
\begin{eqnarray*}
\vert m(t,r) \vert & \leq & \vert m(0,r) \vert + \left \vert
\int_0^t \int_{\vert y \vert \leq r} \int_{\vert v \vert \leq Q_T}
\partial_s (g(s, y, v)) \ dv \ dy \ ds \ \right \vert \\
& = & \vert m(0,r) \vert + \left \vert \int_0^t \int_{\vert y \vert
\leq r} \int_{\vert v \vert \leq Q_T} v \cdot \nabla_y g(s,
y, v) \ dv \ dy \ ds \ \right \vert \\
& \leq & \vert m(0,r) \vert + \int_0^t \int_{\vert y \vert = r}
\int_{\vert v \vert \leq Q_T} \vert v \vert \ \vert g(s,\vert y
\vert, \vert v \vert, y \cdot v) \vert \ dv \ dS_y \ ds \\
& \leq & M(0) + C \int_0^t \int_{\vert y \vert = r} \int_{\vert v
\vert \leq Q_T} \vert v \vert \left (\frac{1 + \int_0^s M(\tau) \ d\tau}{\vert y \vert^2} \right) \ dv \ dS_y  \ ds \\
& = & M(0) + C \int_0^t Q_T^4 (1 + \int_0^s M(\tau) \ d\tau) \ ds \\
& \leq & C + C \int_0^t M(\tau) \ d\tau \\
\end{eqnarray*}
for $ r \geq 2Q_T T$. Then, since $m(t,r) \leq C$ for $ r \leq
2Q_TT$, we find, for all $r$,
$$ m(t,r) \leq C + C \int_0^t M(\tau) \ d\tau$$
and consequently
$$ M(t) \leq C + C \int_0^t M(\tau) \ d\tau.$$
By Gronwall's Inequality,
\begin{equation}
\label{m} \mathcal{M}(t) \leq (\mathcal{M}(0) + CT) (1 + Ct
\exp(Ct)) \leq (\mathcal{M}(0)+ CT) (1 + CT \exp(CT)) \leq C.
\end{equation}
Notice, then, that this bound and (\ref{g}) imply
$$ \vert g(t,x,v) \vert \leq \frac{C}{\vert x \vert^{2}} $$
for $\vert x \vert \geq 2Q_T T$.\\

Proceeding in the standard way, we next estimate the gradient of the
electric field.

\begin{eqnarray*}
\partial_{x_i} E_k & = & \frac{\partial}{\partial x_i} \left (\frac{x_k}{\vert x
\vert^3} \int_{\vert y \vert \leq \vert x \vert} \rho(t,y) dy \right ) \\
& = & \frac{\partial}{\partial x_i}\left (\frac{x_k}{\vert x
\vert^3} \right) m(t,\vert x \vert) + \frac{x_k}{\vert x \vert^3}
\frac{\partial}{\partial x_i} (m(t,\vert x \vert)).
\end{eqnarray*}
Then,
$$ \left\vert \frac{\partial}{\partial x_i} \left (\frac{x_k}{\vert x \vert^3}\right ) \right\vert = \left \vert
\frac{\delta_{ik}}{\vert x \vert^3} - \frac{3x_i x_k}{\vert x
\vert^5} \right \vert \leq \frac{C}{\vert x \vert^3} $$ and,
letting $r = \vert x \vert$,
$$ \left \vert \frac{\partial}{\partial x_i} m(t,r) \right \vert = \left \vert m_r(t,r)
\left (\frac{x_i}{r} \right) \right \vert \leq C \vert \rho(t,r)
\vert r^2. $$ Thus,
\begin{eqnarray*}
\vert \partial_{x_i} E_k \vert & \leq & \vert m(t,r) \vert
\frac{C}{r^3} + \frac{C \vert \rho(t,r) \vert \ r^2}{r^2} \\
& \leq & \frac{C}{r^3} + C \vert \rho(t,r) \vert.
\end{eqnarray*}
Since the best we can do a priori is $\vert \rho(t,r) \vert \leq
Cr^{-2}$, we can only conclude

\begin{equation}
\label{gradfield}
\vert \partial_{x_i} E_k \vert \leq C \vert x
\vert^{-2}.
\end{equation}

We begin to estimate the spatial decay of $\rho$ by first estimating
the large $\vert x \vert$ behavior of the field integral obtained by
integrating (\ref{FI}) in $v$. Assume $\vert x \vert \geq 8Q_TT$ and
$f(t,x,v)$ is nonzero.  Define
$$ \mathcal{E}(t,x) = E(t,x) + A(t,x). $$
Then,
\begin{eqnarray*}
\left \vert \int_{\vert v \vert \leq Q_T} \mathcal{E}(s,X(s))
\cdot \nabla_v F(V(s)) dv \right \vert & \leq & \left \vert
\int_{\vert v \vert \leq Q_T} \mathcal{E}(s,X(s)) \cdot (\nabla_v
F(V(s)) - \nabla_v F(v)) \ dv \right \vert \\
& \ & + \ \left \vert \int_{\vert v \vert \leq Q_T}
(\mathcal{E}(s,X(s)) - \mathcal{E}(s, x + (s-t)v)) \cdot \nabla_v
F(v) \ dv \right \vert \\
& \ & + \ \left \vert \int_{\vert v \vert \leq Q_T} \nabla_v \cdot
(F(v) \ \mathcal{E}(s,x + (s-t)v)) \ dv \right \vert \\
& \ & + \ \left \vert \int_{\vert v \vert \leq Q_T} F(v) \nabla_v
\cdot (\mathcal{E}(s,x + (s-t)v)) \ dv \right \vert \\
& =: & I + II + III + IV.
\end{eqnarray*}

By the Mean Value Theorem, (\ref{e}), (\ref{xchar1}), and (\ref{m}),
\begin{equation}
\label{eq1}
\begin{array}{rcl} I & \leq & \int_{\vert v \vert \leq
Q_T} \vert \mathcal{E}(s,X(s)) \vert \ \Vert \nabla^2_v F
\Vert_\infty \ \vert V(s) - v \vert \ dv \\
\\
& \leq & \int_{\vert v \vert \leq Q_T} (C \vert X(s) \vert^{-2}
\Vert \nabla^2_v F \Vert_\infty \ \vert \int_s^t E(\tau, X(\tau)) \ d\tau \vert ) \ dv \\
\\
& \leq & C Q_T^3 \vert x \vert^{-2} (t-s) \vert x \vert^{-2} \\
\\
& \leq & C \vert x \vert^{-4}.
\end{array}
\end{equation}

To estimate $II$, we again use the Mean Value Theorem and
(\ref{gradfield}), so that there is $\xi_x$ between $X(s)$ and $x +
(s - t)v$ with

\begin{eqnarray*}
\vert \mathcal{E}_i(s, X(s)) - \mathcal{E}_i(s, x + (s - t)v)
\vert  & = & \vert \nabla_x \mathcal{E}_i (s, \xi_x) \cdot (X(s) - x - (s - t)v) \vert \\
& \leq & C \vert \xi_x \vert^{-2} \vert X(s) - x + (t - s)v \vert.
\end{eqnarray*}
Thus, we find
\begin{eqnarray*}
II & \leq & C \int_{\vert v \vert \leq Q_T} \vert \xi_x \vert^{-2}
\ \vert X(s) - (x + (s-t)v) \vert
\ \Vert \nabla_v F \Vert_\infty \ dv \\
& \leq & C \int_{\vert v \vert \leq Q_T} \vert \xi_x \vert^{-2}
\left (\int_s^t \vert V(\tau) - v \vert \ d\tau \right ) \ dv \\
& \leq & C T^2 Q_T^3 \vert x \vert^{-2} \vert \xi_x \vert^{-2},
\end{eqnarray*}
and by (\ref{xchar1}),
\begin{eqnarray*}
\vert \xi_x \vert & \geq & \vert X(s) \vert - \vert X(s) - (x +
(s-t)v) \vert \\
& \geq & \frac{1}{2}\vert x \vert - \left \vert \int_s^t (V(\tau) - v) \ d\tau \right \vert \\
& \geq & \frac{1}{2}\vert x \vert - 2TQ_T \\
& \geq & \frac{1}{4}\vert x \vert.
\end{eqnarray*}
Therefore,
\begin{equation}
\label{eq2}
II \leq C \vert x \vert^{-4}.
\end{equation}
By the Divergence Theorem, we find
\begin{equation}
\label{eq3}
III = 0.
\end{equation}
Then, evaluating $IV$ yields
\begin{eqnarray*}
IV & \leq & \int_{\vert v \vert \leq Q_T} \Vert F \Vert_\infty \
\vert \nabla_x \cdot \mathcal{E}(s,x + (s-t)v) \vert \ \vert
s-t \vert \ dv \\
\\
& \leq & CT \int_{\vert v \vert \leq Q_T} \vert
\rho(s,x + (s-t)v) + \nabla_x \cdot A(s, x + (s-t)v) \vert \ dv \\
& \leq & C \int_{\vert v \vert \leq Q_T} \vert \rho(s,x + (s-t)v)
\vert \ dv + C R^{-4} (x + (s-t)v).
\end{eqnarray*}
Since
$$ \vert x + (s-t)v \vert \geq \vert x \vert - TQ_T \geq \frac{1}{2} \vert x
\vert,$$ we have

\begin{equation}
\label{eq4}
IV \leq C \int_{\vert v \vert \leq Q_T} \vert \rho(s,x
+ (s-t)v) \vert \ dv + C \vert x \vert^{-4}.
\end{equation}

Finally, collecting (\ref{eq1}), (\ref{eq2}), (\ref{eq3}),
(\ref{eq4}), we have for $\vert x \vert \geq 8Q_TT$,
\begin{equation}
\label{eqmain} \left \vert \int_{\vert v \vert \leq Q_T}
\mathcal{E}(s,X(s)) \cdot \nabla_v F(V(s)) \ dv \right \vert \leq C \vert x
\vert^{-4} + C \int_{\vert v \vert \leq Q_T} \vert \rho(s, x +
(s-t)v) \vert \ dv.
\end{equation}

Now, we can bound $\Vert \rho(t) \Vert_4$. Using the bound on the
velocity support, we have

\begin{eqnarray*}
\Vert \rho(t) \Vert_\infty & = & \sup_x \left \vert \int_{\vert v
\vert \leq Q_T} (F(v) - f(t,x,v)) \ dv \right \vert \\
& \leq & \frac{4\pi}{3} Q_T^3 \ \Vert F - f_0 \Vert_\infty \\
& \leq & C.
\end{eqnarray*}
Thus, for $\vert x \vert < 8 Q_T T$, we have
$$ \vert \rho(t,x) \vert \leq \Vert \rho(t) \Vert_\infty \
\left (\frac{8 Q_T T}{\vert x \vert} \right )^4 \leq C \vert x
\vert^{-4}.$$

Now, recall (\ref{FI}) :
$$ g(t,x,v) = g(0,X(0,t,x,v),V(0,t,x,v)) - \int_0^t
\mathcal{E}(s,X(s,t,x,v)) \cdot \nabla_v F(V(s,t,x,v)) \ ds.$$ We
use this equation and (\ref{eqmain}) so that for $\vert x \vert \geq
8Q_T T$,
\begin{eqnarray*}
\vert \rho(t,x) \vert & = & \left \vert \int_{\vert v \vert \leq
Q_T} g(t,x,v) \ dv \right \vert \\
& \leq & \int_{\vert v \vert \leq Q_T} \left ( \ \vert g_0(X(0),
V(0)) \vert + \left \vert \int_0^t \mathcal{E}(s, X(s)) \cdot
\nabla_v F(V(s)) \ ds \right \vert \ \right ) \ dv \\
& \leq & C Q_T^3 \vert x \vert^{-4} + \int_0^t \left ( C Q_T^3
\vert x \vert^{-4} + C \int_{\vert v \vert \leq Q_T} \vert \rho(s,
x + (s-t)v) \vert \ dv \right ) \ ds \\
& \leq & C \vert x \vert^{-4} + C \int_0^t \int_{\vert v \vert
\leq Q_T} \vert \rho(s, x + (s-t)v) \vert \ dv \ ds.
\end{eqnarray*}
So,
$$ \vert x \vert^4 \vert \rho(t,x) \vert \leq C + C \int_0^t \int_{\vert v \vert
\leq Q_T} \vert x \vert^4 \ \vert \rho(s, x + (s-t)v) \vert \ dv \
ds. $$ Define $$ \mathcal{P}(t) := \sup_x ( \vert x \vert^4 \ \vert
\rho(t,x) \vert ). $$ Then, for all $x \in \mathbb{R}^3$,
\begin{eqnarray*}
\mathcal{P}(t) & \leq &  C + C \int_{\vert v \vert \leq Q_T}
\int_0^t \mathcal{P}(s) \ ds \ dv\\
& \leq & C + C \int_0^t \mathcal{P}(s) \ ds
\end{eqnarray*}
and using the Gronwall Inequality, we find
$$ \mathcal{P}(t) \leq C$$
and thus $ \vert \rho(t,x) \vert \leq C \vert x \vert^{-4} $ for all
$x \in \mathbb{R}^3$.  Finally, we may apply Theorem $2$ with any $q
> 7 + \sqrt{33}$, and since $\Vert \rho(t) \Vert_\infty$ is
bounded, we conclude that $\vert \rho(t,x) \vert \leq CR^{-4}(x)$.
Since this estimate is independent of $T$, we find
$$\sup_{t \in [0,T]} \Vert \rho(t) \Vert_p \leq C.$$  This remains
true for any $T > 0$, so the proof of Lemma $1$, and thus Theorem
$1$, is complete.


\section*{Section 2}

As in Section $1$, we will show that the charge density decays at a
faster rate than previously known.  In the work that follows, we
will use the framework of the previous sections and assume
conditions $(I)$ and $(II)$ from Section $1$, without the spherical
symmetry of $f_0$.  However, we will not take condition $(III)$ as
an assumption, and instead, assume condition $(IV)$ holds for some
$C > 0$ and all $t \geq 0$, $x \in \mathbb{R}^3$, and $v \in
\mathbb{R}^3$.  Since we have made a change in the assumptions, we
may not use results from the previous sections, unless otherwise
stated. Thus, this section can be viewed independently from the
others, as we will rely more on results shown previously in
\cite{VPSSA}. Recall from the introduction that we will use $C_{p,
t}$ to denote a
generic constant which depends upon $\Vert \rho(t) \Vert_p$ and $t$.\\

To begin, we apply Theorems $1$ and $2$ of \cite{VPSSA}, finding a
unique $f \in C^1([0,T] \times \mathbb{R}^3 \times \mathbb{R}^3)$
which satisfies (\ref{one}) and assume
\begin{equation}
\label{triple}
\vert \vert \vert (F - f) (t) \vert \vert \vert <
\infty
\end{equation}
for all $t \in [0,T)$. To first bound the velocity support, we
write, as before,
$$ g(t,x,v)  := F(v) - f(t,x,v) $$
and
$$ \mathcal{E}(t,x) := E(t,x) + A(t,x). $$
Using Lemma $1$ of \cite{VPSSA}, we find for all $x \in
\mathbb{R}^3$ and $t \in [0,T)$,
$$ \int_0^t \vert \mathcal{E}(\tau,x) \vert  \ d\tau \leq C \sup_{\tau \in [0,t]}\Vert \rho(\tau)
\Vert_p. $$ Then, for any $s \in [0,t)$,
\begin{eqnarray*}
\vert V(s) - V(0) \vert & \leq & \int_0^s \vert
\mathcal{E}(\tau, X(\tau, t, x, v)) \vert \ d\tau \\
& \leq & C_{p, s}.
\end{eqnarray*}
Assuming $ f \neq 0$, we know $f(t,x,v) = f_0(X(0),V(0))$, and by
$(II)$ it must follow that $$\vert V(0, t, x, v) \vert \leq C.$$
Thus, $\vert V(s) \vert \leq C_{p, s}$ for any $s \in [0,t)$.
Following previous notation let us write
$$Q_t := \sup \{ \vert v \vert : \exists x \in \mathbb{R}^3, \tau \in [0,t) \ \mathrm{such \
that} \ f(\tau, x, v) \neq 0 \} $$ so that for all $s \in [0,t)$

\begin{equation}
\label{velsupp}
\vert V(s) \vert \leq Q_t.
\end{equation}
This bound on the velocity establishes some of the relations shown
in previous sections. Most importantly, (\ref{xchar1}) and
(\ref{xchar2}) must hold for $\vert x \vert \geq 2Q_tt$.\\

Next, we denote the $v$-derivatives of characteristics by
$$ \mathbb{A}_{ij} := \frac{\partial X_i}{\partial v_j}(s,t,x,v) $$
and
$$ \mathbb{B}_{ij} := \frac{\partial V_i}{\partial v_j}(s,t,x,v). $$
Then, using the characteristic equations of (\ref{char}),
\begin{equation}\left. \begin{array}{ccc}
\label{char} & & \frac{\partial \mathbb{A}}{\partial s}(s) = \mathbb{B}(s), \\
 & & \frac{\partial \mathbb{B}}{\partial s}(s) = \nabla_x \mathcal{E}(s, X(s)) \
\mathbb{A}(s), \\
& & \mathbb{A} \vert_{s = t} = 0, \\
& & \mathbb{B} \vert_{s = t} = \mathbb{I}. \end{array} \right\}
\end{equation}
Thus,
$$ \vert \mathbb{A}(s) \vert \leq \int_s^t \vert \mathbb{B}(\tau) \vert \ d\tau $$
and
$$ \vert \mathbb{B}(s) \vert \leq 1 + \int_s^t \vert \nabla_x
\mathcal{E}(\tau, X(\tau)) \vert \ \vert \mathbb{A}(\tau) \vert \
d\tau $$ which leads to
\begin{equation}
\label{AB}
\vert \mathbb{A}(s) \vert + \vert \mathbb{B}(s) \vert
\leq 1 + \int_s^t ( \vert \mathbb{B}(\tau) \vert + \vert \nabla_x
\mathcal{E}(\tau, X(\tau)) \vert \ \vert \mathbb{A}(\tau) \vert) \
d\tau.
\end{equation}
Again, using Lemma $1$ from \cite{VPSSA}, we find
\begin{equation}
\label{DE} \vert \nabla_x \mathcal{E}(\tau, X(\tau)) \vert \leq
C_{p, \tau} R^{-3}(X(\tau)).
\end{equation}
So, define
\begin{equation}
\label{H} \mathcal{H}(s) := \sup_{\{(x,v) : \vert x \vert > 2 Q_t
t \}} ( \vert \mathbb{A}(s,t,x,v) \vert + \vert
\mathbb{B}(s,t,x,v) \vert ).
\end{equation}
Then, for $\vert x \vert > 2Q_tt$, we use (\ref{xchar1}) to find
$$R^{-3}(X(\tau)) \leq R^{-3}(Q_tt),$$ and thus, from (\ref{AB})
$$ \mathcal{H}(s) \leq 1 + \int_s^t \max \{1, C_{p, \tau} R^{-3}(Q_t t) \} \ \mathcal{H}(\tau) \ d\tau. $$
Using the Gronwall Inequality, we find
\begin{equation}
\label{H2}
\mathcal{H}(s) \leq C_{p, t}
\end{equation}
for $s \in [0,t]$, and v-derivatives of both characteristics are bounded
for $\vert x \vert > 2Q_tt$.\\

Summarizing (\ref{char}), we may write
$$ \frac{\partial^2 \mathbb{A}}{ds^2}(s)  = \nabla_x \mathcal{E}(s,X(s)) \mathbb{A}(s); \ \ \
\mathbb{A} \vert_{s = t} = 0; \ \ \ \frac{\partial
\mathbb{A}}{\partial s} \vert_{s = t} = \mathbb{I}.$$ So,
\begin{eqnarray*}
\mathbb{A}(s) & = & (s-t) \mathbb{I} + \int_s^t \int_\tau^t
\nabla_x \mathcal{E}(\lambda, X(\lambda)) \mathbb{A}(\lambda) \
d\lambda d\tau \\
& =: & (s-t)\mathbb{I} + \gamma_1(s,t,x,v).
\end{eqnarray*}
Again using (\ref{char}),
\begin{eqnarray*}
\frac{\partial \mathbb{B}}{\partial s}(s) & = & \nabla_x \mathcal{E}(s, X(s)) \mathbb{A}(s) \\
& = & \nabla_x \mathcal{E}(s,X(s)) ((s-t)\mathbb{I} + \gamma_1(s,t,x,v)) \\
& =: & (s-t) \nabla_x \mathcal{E}(s,X(s)) + \gamma_2(s,t,x,v).
\end{eqnarray*}
Thus, $$ \mathbb{B}(s) = \mathbb{I} + \gamma_3(s,t,x,v)
$$ where $$\gamma_3(s,t,x,v) = - \int_s^t ((\tau - t) \nabla_x
\mathcal{E}(\tau,X(\tau)) + \gamma_2(\tau,t,x,v) ) \ d\tau.$$ Now,
\begin{eqnarray*}
\frac{\partial}{\partial s} (\det \mathbb{B}) & = & \det \left (
\begin{array}{ccc} \dot{\mathbb{B}}_{11} & \dot{\mathbb{B}}_{12} & \dot{\mathbb{B}}_{13} \\
\mathbb{B}_{21} & \mathbb{B}_{22} & \mathbb{B}_{23} \\
\mathbb{B}_{31} & \mathbb{B}_{32} & \mathbb{B}_{33}
\end{array} \right ) + \det \left (
\begin{array}{ccc} \mathbb{B}_{11} & \mathbb{B}_{12} & \mathbb{B}_{13} \\
\dot{\mathbb{B}}_{21} & \dot{\mathbb{B}}_{22} & \dot{\mathbb{B}}_{23} \\
\mathbb{B}_{31} & \mathbb{B}_{32} & \mathbb{B}_{33}
\end{array} \right ) + \det \left (
\begin{array}{ccc} \mathbb{B}_{11} & \mathbb{B}_{12} & \mathbb{B}_{13} \\
\mathbb{B}_{21} & \mathbb{B}_{22} & \mathbb{B}_{23} \\
\dot{\mathbb{B}}_{31} & \dot{\mathbb{B}}_{32} &
\dot{\mathbb{B}}_{33}
\end{array} \right ) \\
& =: & I + II + III.
\end{eqnarray*}
Estimating the first term yields
\begin{eqnarray*}
I & = & \det \left ( \begin{array}{ccc} (s-t) \frac{\partial
\mathcal{E}_1}{\partial x_1} + (\gamma_2)_{11} & (s-t)
\frac{\partial \mathcal{E}_1}{\partial x_2} + (\gamma_2)_{12} &
(s-t) \frac{\partial \mathcal{E}_1}{\partial x_3} + (\gamma_2)_{13} \\
(\gamma_3)_{21} & 1 + (\gamma_3)_{22} & (\gamma_3)_{23} \\
(\gamma_3)_{31} & (\gamma_3)_{32} & 1 + (\gamma_3)_{33}
\end{array} \right ) \\
& =: & (s-t) \frac{\partial \mathcal{E}_1}{\partial x_1} +
(\gamma_2)_{11} + \sigma_1(s,t,x,v).
\end{eqnarray*}
Similarly,
$$ II =: (s-t) \frac{\partial \mathcal{E}_2}{\partial x_2} + (\gamma_2)_{22} +
\sigma_2(s,t,x,v) $$ and
$$ III =: (s-t) \frac{\partial \mathcal{E}_3}{\partial x_3} + (\gamma_2)_{33} +
\sigma_3(s,t,x,v). $$ So, we find
$$ \frac{\partial}{\partial s} (\det \mathbb{B}) = (s-t) \nabla_x \cdot
\mathcal{E}(s,x) \vert_{x = X(s)} + \sum_{j=1}^3 ( (\gamma_2)_{jj}
+ \sigma_j) $$
and since $\mathbb{B} \vert_{s = t} = \mathbb{I}$,
$$ \det \mathbb{B} = 1 - 4\pi \int_s^t (\tau - t) \rho(\tau, X(\tau)) \ d\tau
- \int_s^t \sum_{j=1}^3 ( (\gamma_2)_{jj} + \sigma_j) \ d\tau. $$

Let $$ \epsilon(s,t,x,v) := 4\pi \int_s^t (\tau - t) \rho(\tau,
X(\tau)) \ d\tau + \int_s^t \sum_{j=1}^3 ( (\gamma_2)_{jj} +
\sigma_j) \ d\tau. $$ Then, for $\vert \epsilon \vert < 1$,
\begin{equation}
\label{detB}
\begin{array}{rcl}
\frac{1}{\det \mathbb{B}} & = & \frac{1}{1 - \epsilon} \\
\\
& = & 1 + \epsilon + \sum_{n=2}^\infty \epsilon^n \\
\\
& =: & 1 + 4\pi \int_s^t (\tau - t) \rho(\tau, X(\tau)) \ d\tau +
\eta(s,t,x,v).
\end{array}
\end{equation}

Now that the determinant has been written in a nicer form, we
estimate the remaining terms.  Let $\vert x \vert
> 2Q_tt$. Using (\ref{xchar1}), (\ref{triple}), (\ref{DE}),
(\ref{H}), and (\ref{H2}), we find the following bounds on the error
terms for any $i,j = 1, 2, 3$:
\begin{eqnarray*}
\vert (\gamma_1)_{ij} \vert & \leq & C_{p, t} T^2 R^{-3}(x)
\sup_{s \in [0,t]} \vert \mathbb{A}(s,t,x,v) \vert \\
& \leq & C_{p, t} R^{-3}(x),
\end{eqnarray*}

$$ \vert (\gamma_2)_{ij} \vert \leq C_{p, t} R^{-6}(x), $$
and

\begin{eqnarray*}
\vert (\gamma_3)_{ij} \vert & \leq & C_{p, t} T^2 R^{-3}(x) + C_{p, t} R^{-6}(x) \\
& \leq & C_{p, t} R^{-3}(x).
\end{eqnarray*}
Then, using the estimates of the $\gamma$ terms, for any $k =
1,2,3$,
\begin{eqnarray*}
\vert \sigma_k \vert & \leq & C (\vert \nabla_x \mathcal{E} \vert
+ \vert \gamma_2 \vert) (2\vert \gamma_3 \vert + 2 \vert \gamma_3
\vert ^2) + 2(\vert \nabla_x \mathcal{E} \vert + \vert \gamma_2
\vert) (\vert \gamma_3 \vert + 2 \vert \gamma_3 \vert^2)\\
\\
& \leq & [C_{p, t} T R^{-3}(x) + C_{p, t} R^{-6}(x)][2C_{p, t} R^{-3}(x) + 2C_{p, t} R^{-6}(x)] \\
\\
& \ &  \ \ + \ 2[C_{p, t} T R^{-3}(x) +
C_{p, t} R^{-6}][C_{p, t} R^{-3}(x) + 2C_{p, t} R^{-6}(x)] \\
\\
& \leq & C_{p, t} R^{-6}(x).
\end{eqnarray*}

Now, we may bound $\epsilon$.
\begin{eqnarray*}
\vert \epsilon \vert & \leq & 4\pi T^2 \sup_{\tau \in [0,t]} \Vert
\rho(\tau) \Vert_p R^{-p}(x) + 3T \sup_{\tau \in [s,t]}(\vert \gamma_2(\tau) \vert + \vert \sigma(\tau) \vert) \\
& \leq & C_{p, t}  R^{-p}(x) + C_{p, t} R^{-6}(x) \\
& \leq & C_{p, t} R^{-p}(x) \\
& =: & C^{(1)} R^{-p}(x) \\
& < & \frac{1}{2}
\end{eqnarray*}
for $\vert x \vert > (2C^{(1)})^{\frac{1}{p}}$. Thus
$\sum_{n=2}^\infty \epsilon^n$ converges and $\eta$ is well-defined
for large enough values of $\vert x \vert$. Finally, we have, for
$\vert x \vert > \max\{ 2Q_tt, (2C^{(1)})^\frac{1}{p}\}$,
\begin{equation}
\label{eta}
\begin{array}{rcl}
\vert \eta \vert & = & \left \vert \int_s^t \sum_{j=1}^3
((\gamma_2)_{jj} + \sigma_j) \ d\tau + \sum_{n=2}^\infty \epsilon^n \right \vert \\
\\
& \leq & 3C_{p, t} T R^{-6}(x) + ((C^{(1)})^2 R^{-2p}(x) + (C^{(1)})^3 R^{-3p}(x) + ...) \\
\\
& \leq & C_{p, t} R^{-6}(x).
\end{array}
\end{equation}

Using (\ref{velsupp}), we find for $t \in [0,T)$,
\begin{eqnarray*}
\Vert \rho(t) \Vert_\infty & = & \sup_x \left \vert \int (F(v) -
f(t,x,v)) \ dv \right \vert \\
& \leq & \frac{4\pi}{3}Q_t^3 \ \Vert F - f_0 \Vert_\infty \\
& \leq & C_{p,t}.
\end{eqnarray*}
Thus, if $\vert x \vert \leq D$ for some $D > 0$, we have
\begin{equation}
\label{smallx}
\begin{array}{rcl}
\vert \rho(t,x) \vert & \leq & \Vert \rho(t) \Vert_\infty R^6(D) R^{-6}(x) \\
\\
& \leq & C_{p,t} R^{-6}(x).
\end{array}
\end{equation}
Now, denote $C^{(2)} :=  \max \{2Q_t t, (2C^{(1)})^\frac{1}{p}, 2N
\}$, and let $\vert x \vert > C^{(2)} $. Then, by (\ref{xchar1}) we
have
$$ \vert X(0,t,x,v) \vert \geq \frac{1}{2} \vert x \vert > N. $$
Using $(V)$, (\ref{fchar}), (\ref{detB}), and (\ref{eta}) to
estimate $\rho(t,x)$ yields
\begin{eqnarray*}
\int f(t,x,v) \ dv & = & \int f_0(X(0,t,x,v), V(0,t,x,v)) \ dv \\
& = & \int F(V(0,t,x,v)) \ dv \\
& = & \int F(V(0,t,x,v)) (\det \mathbb{B}(0)) \frac{1}{\det \mathbb{B}(0)} \ dv \\
& = & \int F(V(0,t,x,v)) ( 1 + 4\pi \int_0^t (\tau - t) \rho(\tau,
X(\tau,t,x,v)) d\tau + \eta(0,t,x,v) ) (\det \frac{\partial
V}{\partial v} (0,t,x,v) ) \ dv \\
& = & \int F(w) \ dw + \int F(V(0,t,x,v)) \eta(0,t,x,v) (\det
\frac{\partial V}{\partial v} (0,t,x,v) ) \ dv \\
& & \ + 4\pi \int F(V(0,t,x,v)) \int_0^t (\tau - t) \rho(\tau,
X(\tau,t,x,v)) \ d\tau \ (\det \frac{\partial V}{\partial v}
(0,t,x,v) ) \ dv.
\end{eqnarray*}
Since $\rho(t,x) = \int (F(v) - f(t,x,v)) \ dv$, we have

\begin{equation} \begin{array}{rcl}
\label{rho} \rho(t,x) & = & 4\pi \int F(V(0,t,x,v)) \int_0^t (t -
\tau) \rho(\tau, X(\tau,t,x,v)) \ d\tau \ ( \det \frac{\partial
V}{\partial v} (0,t,x,v) ) \ dv \\
\\
& & \ - \int F(V(0,t,x,v)) \eta(0,t,x,v) (\det \frac{\partial
V}{\partial v} (0,t,x,v) ) \ dv.
\end{array}
\end{equation}
Now, let
$$ \Psi(t) := \Vert \rho(t) \Vert_6 = \sup_x ( \vert \rho(t,x) \vert R^6(x)). $$
In order to make the change of variables $w = V(0,t,x,v)$ in the
$v$-integral, we must first show that the mapping $v \rightarrow
V(0,t,x,v)$ is bijective. For the moment, we will take this for
granted, and continue with the estimate of $\rho$, delaying the
proof of this fact until the very end.  By (\ref{eta}), (\ref{rho}),
and the work of Section $3.6$, we have for $\vert x \vert >
C^{(2)}$,
\begin{equation} \begin{array}{rcl}
\label{rho2} \vert \rho(t,x) \vert & \leq & 4\pi \int
F(V(0,t,x,v)) (\int_0^t (t - \tau) \Psi(\tau)
R^{-6}(X(\tau,t,x,v)) \ d\tau)
(\det \frac{\partial V}{\partial v} (0,t,x,v) ) \ dv \\
\\
& & \ + C_{p, t} R^{-6}(x) \int F(w) \ dw \\
\\
& \leq & CR^{-6}(x) \int F(V(0,t,x,v)) (\det \frac{\partial
V}{\partial v} (0,t,x,v) ) dv \int_0^t (t - \tau) \Psi(\tau) \ d\tau  \\
\\
& & \ + C_{p, t} R^{-6}(x) \\
\\
& \leq & R^{-6}(x) ( C \int_0^t (t - \tau) \Psi(\tau) \ d\tau +
C_{p, t} ).
\end{array}
\end{equation}
So, applying (\ref{smallx}) with $D = C^{(2)}$ and combining with
(\ref{rho2}), we have for all $x$,
$$ \vert \rho(t,x) \vert R^6(x) \leq C_{p, t} + C \int_0^t (t-\tau) \Psi(\tau) d\tau $$ and
since the right side is independent of $\vert x \vert$,
$$ \Psi(t) \leq C_{p, t} + C \int_0^t (t-\tau) \Psi(\tau) d\tau.$$
By the Gronwall Inequality,
$$ \Psi(t) \leq C_{p, t}.$$
Therefore, for every $t \in [0,T)$,
$$ \Vert \rho(t) \Vert_6 \leq C_{p, t}$$ and the proof is complete.

\section*{Change of Variables}

In order to justify the change of variables used in line
(\ref{rho2}) of the previous section, we must first demonstrate that
the mapping $v \rightarrow V(0,t,x,v)$ is bijective. This is done
below.\\

From (\ref{DE}), we know there is $C^{(3)}_{p,t} > 0$ such that
$$\vert \nabla \mathcal{E}(t,x) \vert \leq C^{(3)}_{p,t} \vert x \vert^{-3}$$
for every $x \in \mathbb{R}^3$, $t \in [0,T]$.  Also, using
(\ref{H2}), we know there is $C^{(4)}_{p,t} > 0$ such that for all
$s,t \in [0,T]$, $x,v \in \mathbb{R}^3$, with $\vert x \vert > 2T
Q(T)$,
$$ \left \vert \frac{\partial X}{\partial v}(s,t,x,v) \right \vert \leq C^{(4)}_{p,t}.$$
Finally, using Lemma $1$ of \cite{VPSSA}, we know there is
$C^{(5)}_{p,t} > 0$ such that for all $x \in \mathbb{R}^3$ and $t
\in [0,T]$,
$$ \int_0^t \vert \mathcal{E} (\tau,x) \vert \ d\tau \leq C^{(5)}_{p,t}. $$
For any $D > 0$, define $C_D := \max \{8 \left( D + C^{(5)}_{p,t} T
\right), 2TQ(T), 4 (6 C^{(4)}_{p,t} C^{(3)}_{p,t} T)^\frac{1}{3}\}$
and $B(0,D) : = \{v \in \mathbb{R}^3 : \vert v \vert \leq D\}$.
Injectivity on $B(0,D)$ can now be shown in the following lemma.

\begin{lemma}
For any $D > 0 $, $\vert x \vert > C_D$, and $t \in [0,T]$, the
mapping $v \rightarrow V(0,t,x,v)$ is injective on $B(0,D)$.
\end{lemma}

\noindent {\bf Proof:} Let $D > 0$, $\vert x \vert > C_D$ and $t \in
[0,T]$ be given. Then, let $v_1, v_2 \in B(0,D)$ be given with
$$ V(0,t,x,v_1) = V(0,t,x,v_2).$$
We have  $$ v_1 + \int_0^t \mathcal{E}(\tau, X(\tau, t, x, v_1)) \
d\tau = v_2 + \int_0^t \mathcal{E}(\tau, X(\tau, t, x, v_2)) \
d\tau. $$ So, using the Mean Value Theorem, for any $i = 1,2,3$
there is $\xi_x^i$ between $X(\tau,t,x,v_1)$ and $X(\tau,t,x,v_2)$
and $\xi_v^i$ between $v_1$ and $v_2$ such that

$$ \mathcal{E}_i(\tau, X(\tau,t,x,v_1)) - \mathcal{E}_i(\tau, X(\tau,t,x,v_2)) =
\nabla_x \mathcal{E}_i(\tau,\xi_x^i) \cdot \left ( X(\tau,t,x,v_1)
- X(\tau,t,x,v_2) \right )$$ and
$$ X_i(\tau,t,x,v_1) - X_i(\tau,t,x,v_2) = \nabla_v
X_i(\tau,t,x,\xi_v^i) \cdot (v_2 - v_1).$$ Then, since
\begin{eqnarray*}
\vert \xi_x^i \vert & \geq & \vert X(\tau,t,x,v_1) \vert - \vert X(\tau,t,x,v_1) - X(\tau,t,x,v_2) \vert \\
& \geq & \frac{1}{2} \vert x \vert - \int_\tau^t \vert V(\iota,t,x,v_1) - V(\iota,t,x,v_2) \vert \ d\iota \\
& \geq & \frac{1}{2} \vert x \vert - \int_\tau^t \left ( \vert
V(\iota,t,x,v_1) - v_1 \vert + \vert v_1 - v_2 \vert
+ \vert V(\iota,t,x,v_2) - v_2 \vert \right ) \ d\iota \\
& \geq & \frac{1}{2} \vert x \vert - \left ( \vert v_1 - v_2 \vert
+ 2 \int_\tau^t \int_0^\iota \Vert E(\lambda) \Vert_\infty \ d\lambda \ d\iota \right ) \\
& \geq & \frac{1}{2} \vert x \vert - 2 \left ( D + C^{(5)}_{p,t}T \right ) \\
& \geq & \frac{1}{4} \vert x \vert,
\end{eqnarray*}
and

$$ \vert X_i(\tau,t,x,v_1) - X_i(\tau,t,x,v_2) \vert \leq C^{(4)}_{p,t} \ \vert v_2 - v_1
\vert,$$ for any $i=1,2,3$, we find

\begin{eqnarray*}
 \vert \mathcal{E}_i(\tau, X(\tau,t,x,v_1)) - \mathcal{E}_i(\tau,
X(\tau,t,x,v_2))\vert & \leq & \vert \nabla_x \mathcal{E}_i(\tau,
\xi_x^i) \vert \ \vert X(\tau,t,x,v_1) - X(\tau,t,x,v_2) \vert \\
& \leq & \left ( \sup_i C^{(3)}_{p,t} \vert \xi_x^i \vert^{-3} \right ) \ \left ( \sqrt{3} C^{(4)}_{p,t} \  \ \vert v_1 - v_2 \vert \right ) \\
& \leq & \sqrt{3} C^{(3)}_{p,t} C^{(4)}_{p,t} \left ( \frac{1}{4}
\vert x \vert \right)^{-3} \vert v_1 - v_2 \vert.
\end{eqnarray*}
Therefore, we have
\begin{eqnarray*}
\vert v_1 - v_2 \vert & \leq & \int_0^t \vert \mathcal{E}(\tau,
X(\tau,t,x,v_1)) - \mathcal{E}(\tau, X(\tau,t,x,v_2)) \vert \ d\tau \\
& \leq & 3 C^{(4)}_{p,t} \vert v_1 - v_2 \vert \  C^{(3)}_{p,t}
\int_0^t \left ( \frac{1}{4} \vert x \vert \right)^{-3} \ d\tau. \\
& \leq & (3 C^{(3)}_{p,t} C^{(4)}_{p,t} T (\frac{1}{4} \vert x \vert)^{-3}) \vert v_1 - v_2 \vert \\
& \leq & \frac{1}{2} \vert v_1 - v_2 \vert.
\end{eqnarray*}
Thus, $\vert v_1 - v_2 \vert = 0$, which implies $v_1 = v_2$, and injectivity follows.\\

\begin{flushright}
$\Box$
\end{flushright}

Now, let $t \in [0,T]$ and $x \in \mathbb{R}^3$ be given.  Define
$$S := \{ w : F(w) \neq 0 \}$$ and $$V^{-1}(S) := \{v : F(V(0,t,x,v))
\neq 0 \}.$$ Using the compact support of $F$, we conclude that $v
\in V^{-1}(S)$ implies
$$ \vert v \vert  \leq W + \int_0^t \Vert \mathcal{E}(\tau) \Vert_\infty \ d\tau \leq C_{p,t}. $$
Therefore, there is a $D > 0$, such that $V^{-1}(S) \subset B(0,D)$.
Thus, for $\vert x \vert > C_D$, and $t \in [0,T]$, the mapping $v
\rightarrow V(0,t,x,v)$ is injective on $V^{-1}(S)$ and bijective
from $V^{-1}(S)$ to $S$. Finally,

\begin{eqnarray*}
\int F(V(0,t,x,v)) \det \left (\frac{\partial V}{\partial v} \right)
\ dv & = &
\int_{V^{-1}(S)}  F(V(0,t,x,v)) \det \left (\frac{\partial V}{\partial v} \right ) \ dv \\
& = & \int_S F(w) \ dw \\
& = & \int F(w) \ dw.
\end{eqnarray*}
Thus, the change of variables is valid and the justification is
complete.

\end{document}